\documentclass[a4paper,11pt]{article}
\usepackage[latin1]{inputenc}
\usepackage[T1]{fontenc}
\usepackage[english]{babel}
\usepackage{amssymb}
\usepackage{amsfonts, amsmath}
\newtheorem{Theo}{Theorem}

\newtheorem{Rem}{Remark}[section]

\title{Prym varieties and applications}
\author{\textbf{A. Lesfari}
\\\emph{Department of Mathematics}
\\\emph{Faculty of Sciences}
\\\emph{University of Choua\"{i}b Doukkali}
\\\emph{B.P. 20, El-Jadida, Morocco}.
\\\emph{E. mail address} : Lesfariahmed@yahoo.fr, lesfari@ucd.ac.ma}
\date{}
\begin{document}
\maketitle
\begin{abstract}
The classical definition of Prym varieties deals with the
unramified covers of curves. The aim of the present paper is to
give explicit algebraic descriptions of the Prym varieties
associated to ramified double covers of algebraic curves. We make
a careful study of the connection with the concept of algebraic
completely integrable systems and we apply the methods to some
problems of Mathematical Physics.\\
\emph{2000 MSC}: 14Q05,14H40,70H06.\\
\emph{Subj. Class.}: algebraic geometry, Classical integrable systems.\\
\emph{Keywords}: Curves, Prym varieties, integrable systems.\\
\end{abstract}

\vskip0.4cm

\section{Introduction}

During the last decades, algebraic geometry has become a tool for
solving differential equations and spectral questions of mechanics
and mathematical physics. The present paper consists of two
separate but related topics : the first part purely
algebraic-geometric, the second one on Prym varieties in algebraic
integrability.

Prym variety $\Pr ym(\mathcal{C}/\mathcal{C}_{0})$\ \ is a
subabelian variety of the jacobian variety
\begin{eqnarray}
Jac(\mathcal{C})&=&Pic^{0}(\mathcal{C}),\nonumber\\
&=&H^{1}(\mathcal{O}_{C})/ H^{1}(\mathcal{C},\mathbb{Z}),\nonumber
\end{eqnarray}
constructed from a double cover $\mathcal{C}$ of a curve
$\mathcal{C}_0$: if $\sigma$ is the involution on $C$
interchanging sheets, then $\sigma$ extends by linearity to a map
$$\sigma :Jac(\mathcal{C})\longrightarrow Jac(\mathcal{C}).$$ Up to
some points of order two, $Jac(\mathcal{C})$ splits into an even
part and an odd part : the even part is $Jac(\mathcal{C}_0)$ and
the odd part is a $Prym(\mathcal{C}/\mathcal{C}_{0})$. The
classical definition of Prym varieties deals with the unramified
double covering of curves and was introduced by W. Schottky and H.
W. E. Jung in relation with the Schottky problem [6] of
characterizing jacobian varieties among all principally polarized
abelian varieties (an abelian variety is a complex torus that can
be embedded into projective space). The theory of Prym varieties
was dormant for a long time, until revived by D. Mumford around
1970. It now plays a substantial role in some contemporary
theories, for example integrable systems [1,2,5,8,10,14,15,17],
the Kadomtsev-Petviashvili equation( KP equation), in the
deformation theory of two-dimensional Schrödinger operators [19],
in relation to Calabi-Yau three-folds and string theory, in the
study of the generalized theta divisors on the moduli spaces of
stable vector bundles over an algebraic curve [7,11],...

Integrable hamiltonian systems are nonlinear ordinary differential
equations described by a hamiltonian function and possessing
sufficiently many independent constants of motion in involution.
By the Arnold-Liouville theorem [4,16], the regular compact level
manifolds defined by the intersection of the constants of motion
are diffeomorphic to a real torus on which the motion is
quasi-periodic as a consequence of the following differential
geometric fact; a compact and connected $n$-dimensional manifold
on which there exist $n$ vector fields which commute and are
independent at every point is diffeomorphic to an $n$-dimensional
real torus and each vector field will define a linear flow there.
A dynamical system is algebraic completely integrable (in the
sense of Adler-van Moerbeke [1]) if it can be linearized on a
complex algebraic torus $\mathbb{C}^{n}/lattice$ (=abelian
variety). The invariants (often called first integrals or
constants) of the motion are polynomials and the phase space
coordinates (or some algebraic functions of these) restricted to a
complex invariant variety defined by putting these invariants
equals to generic constants, are meromorphic functions on an
abelian variety. Moreover, in the coordinates of this abelian
variety, the flows (run with complex time) generated by the
constants of the motion are straight lines.

One of the remarkable developments in recent mathematics is the
interplay between algebraic completely integrable systems and Prym
varieties. The period of these Prym varieties provide the exact
periods of the motion in terms of explicit abelian integrals. The
aim of the first part of the present paper is to give explicit
algebraic descriptions of the Prym varieties associated to
ramified double covers of algebraic curves. The basic algebraic
tools are known and can be found in the book by Arbarello,
Cornalba, Griffiths, Harris [3] and in Mumford's paper [18]. In
the second part of the paper, we make a careful study of the
connection with the concept of algebraic completely integrable
systems and we apply the methods to some problems such as the
Hénon-Heiles system, the Kowalewski rigid body motion and
Kirchhoff's equations of motion of a solid in an ideal fluid. The
motivation was the excellent Haine's paper [10] on the integration
of the Euler-Arnold equations associated to a class of geodesic
flow on $SO(4)$. The Kowalewski's top and the Clebsch's case of
Kirchhoff's equations describing the motion of a solid body in a
perfect fluid were integrated in terms of genus two hyperelliptic
functions by Kowalewski [13] and Kötter [12] as a result of
complicated and mysterious computations. The concept of algebraic
complete integrability (Adler-van Moerbeke) throws a completely
new light on these two systems. Namely, in both cases (see [14]
for Kowalewski's top and [10] for geodesic flow on $SO(4)$), the
affine varieties in $\mathbb{C}^6$ obtained by intersecting the
four polynomial invariants of the flow are affine parts of Prym
varieties of genus 3 curves which are double cover of elliptic
curves. Such Prym varieties are not principally polarized and so
they are not isomorphic but only isogenous to Jacobi varieties of
genus two hyperelliptic curves. This was a total surprise as it
was generically believed that only jacobians would appear as
invariants manifolds of such systems. The method that is used for
reveiling the Pryms is due to L.Haine [10]. By now many algebraic
completely integrable systems are known to linearize on Prym
varieties.

\section{Prym varieties}

Let $$\varphi :\mathcal{C}\longrightarrow \mathcal{C}_{0},$$ be a
double covering with 2n branch points where $\mathcal{C}$ and
$\mathcal{C}_0$ are nonsingular complete curves. Let
$$\sigma:\mathcal{C}\longrightarrow \mathcal{C},$$ be the involution
exchanging sheets of $\mathcal{C}$ over
$\mathcal{C}_0=\mathcal{C}/\sigma.$ By Hurwitz's formula, the
genus of $\mathcal{C}$ is $g=2g_0+n-1,$ where $g_0$ is the genus
of $\mathcal{C}_0.$ Let $\mathcal{O}_{\mathcal{C}}$ be the sheaf
of holomorphic functions on $\mathcal{C}$. Let $S^{g}(
\mathcal{C})$ be the g-th symmetric power of $\mathcal{C}$ (the
totality of unordered sets of points of $\mathcal{C}$). On $S^{g}(
\mathcal{C})$, two divisors $\mathcal{D}_{1}$ and
$\mathcal{D}_{2}$ are linearly equivalent (in short,
$\mathcal{D}_{1}\equiv \mathcal{D}_{2}$) if their difference
$\mathcal{D}_{1}-\mathcal{D}_{2}$ is the divisor of a meromorphic
function or equivalently if and only if
$$\int_{\mathcal{D}_{1}}^{\mathcal{D}_{2}}\omega =\int_{\gamma
}\omega ,\text{ }\forall \omega \in \Omega _{\mathcal{C}}^{1},$$
where $\Omega _{\mathcal{C}}^{1}$ is the sheaf of holomorphic
1-forms on $\mathcal{C}$ and $\gamma $ is a closed path on
$\mathcal{C}.$ We define the jacobian (Jacobi variety) of
$\mathcal{C},$ to be $$Jac(\mathcal{C})=S^{g}(\mathcal{C}) /\equiv
.$$ To be precise, the jacobian $Jac(\mathcal{C})=Pic^{0}(
\mathcal{C})$ of $\mathcal{C}$ is the connected component of its
Picard group parametrizing degree 0 invertible sheaves. It is a
compact commutative algebraic group, i.e., a complex torus.
Indeed, from the fundamental exponential sequence, we get an
isomorphism
\begin{eqnarray}
Jac(\mathcal{C})&=&Pic^{0}( \mathcal{C}),\nonumber\\
&\simeq& H^{0}(\mathcal{C},\Omega _{\mathcal{C}}^{1})
^{*}/H_{1}(\mathcal{C},\mathbb{Z}),\nonumber\\
&\simeq& \mathbb{C}^{g}/\mathbb{Z}^{2g},\nonumber
\end{eqnarray}
via the duality given by Abel's theorem. Consider the following
mapping
$$S^{g}\left( \mathcal{C}\right) \longrightarrow
\mathbb{C}^{g}/L_{\Omega },\text{ }\sum_{k=1}^{g}\mu
_{k}\longmapsto \sum_{k=1}^{g}\int^{\mu _{k}\left( t\right)
}\left( \omega _{1},\ldots ,\omega _{g}\right) =t\left(
k_{1},\ldots ,k_{g}\right) ,$$ where $\left( \omega _{1},\ldots
,\omega _{g}\right) $ is a basis of $\Omega _{\mathcal{C}}^{1}$,
$L_{\Omega }$ is the lattice associated to the period matrix
$\Omega $ and $\mu _{1},\ldots ,\mu _{g}$ some appropriate
variables defined on a non empty Zariski open set. Let $(
a_{1},...,a_{g_{0}},...,a_{g},b_{1},...,b_{g_{0}},...,b_{g}),$ be
a canonical homology basis of $H_1(\mathcal{C},\mathbb{Z})$ such
that
$$\sigma \left( a_{1}\right)=a_{g_{0}+n},...,\sigma \left(
a_{g_{0}}\right) =a_{g},$$ $$\sigma
\left(a_{g_{0}+1}\right)=-a_{g_{0}+1},...,\sigma
\left(a_{g_{0}+n-1}\right) =-a_{g_{0}+n-1},$$ $$\sigma \left(
b_{1}\right)=b_{g_{0}+n},...,\sigma \left( b_{g_{0}}\right)
=b_{g},$$ $$\sigma \left(
b_{g_{0}+1}\right)=-b_{g_{0}+1},...,\sigma \left(
b_{g_{0}+n-1}\right)=-b_{g_{0}+n-1},$$ for the involution
$\sigma.$ Notice that $\varphi \left( a_{g_{0}+1}\right) ,...
,\varphi \left( a_{g_{0}+n-1}\right), \varphi \left(
b_{g_{0}+1}\right) ,... ,\varphi \left( b_{g_{0}+n-1}\right) $ are
homologous to zero on $\mathcal{C}_{0}$. Let $\left( \omega
_{1},\ldots ,\omega _{g}\right) $ be a basis of holomorphic
differentials on $\mathcal{C}$ where $\omega _{g_{0}+n},\ldots
,\omega _{g},$ are holomorphic differentials on $\mathcal{C}_{0}$
and
$$\sigma ^{*}(\omega _{j}) = \left\{\begin{array}{rl}
-\omega _{j},&1\leq j\leq g_{0}+n-1,\\
\omega _{j},&g_{0}+n\leq j\leq g,
\end{array}\right.$$
the pullback of $\omega _{j}.$ The period matrix $\Omega $ of
$Jac(\mathcal{C})$ is explicitly given by
$$\Omega=\left(\begin{array}{cccccc}
A&B&C&D&E&F\\
G&H&I&J&K&L
\end{array}\right),
$$ where $A,\ldots ,L$ denote the following matrices :
\begin{equation}\label{eqn:euler}
A=\left(\begin{array}{ccc} \int_{a_{1}}\omega _{1}&\ldots
&\int_{a_{g_{0}}}\omega _{1}\\
\vdots & &\vdots \\
\int_{a_{1}}\omega _{_{g_{0}+n-1}}&\ldots &\int_{a_{g_{0}}}\omega
_{g_{0}+n-1}
\end{array}\right),
\end{equation}
$$B=\left(\begin{array}{ccc} \int_{a_{g_{0}+1}}\omega _{1}&\ldots
&\int_{a_{g_{0}+n-1}}\omega _{1}\\
\vdots & &\vdots \\
\int_{a_{g_{0}+1}}\omega _{g_{0}+n-1}&\ldots
&\int_{a_{g_{0}+n-1}}\omega _{g_{0}+n-1}
\end{array}\right),$$
$$C=\left(\begin{array}{ccc} \int_{a_{g_{0}+n}}\omega _{1}&\ldots
&\int_{a_{g}}\omega _{1}\\
\vdots & &\vdots \\
\int_{a_{g_{0}+n}}\omega _{g_{0}+n-1}&\ldots &\int_{a_{g}}\omega
_{g_{0}+n-1}
\end{array}\right),$$
$$D=\left(\begin{array}{ccc} \int_{b_{1}}\omega _{1}&\ldots
&\int_{b_{g_{0}}}\omega _{1}\\
\vdots & &\vdots \\
\int_{b_{1}}\omega _{_{g_{0}+n-1}}&\ldots &\int_{b_{g_{0}}}\omega
_{g_{0+n-1}}
\end{array}\right),$$
$$E=\left(\begin{array}{ccc} \int_{b_{g_{0}+1}}\omega _{1}&\ldots
&\int_{b_{g_{0}+n-1}}\omega _{1}\\
\vdots & &\vdots \\
\int_{b_{g_{0}+1}}\omega _{_{g_{0}+n-1}}&\ldots
&\int_{b_{g_{0}+n-1}}\omega _{g_{0+n-1}}
\end{array}\right),$$
$$F=\left(\begin{array}{ccc} \int_{b_{g_{0}+n}}\omega _{1}&\ldots
&\int_{b_{g}}\omega _{1}\\
\vdots & &\vdots \\
\int_{b_{g_{0}+n}}\omega _{_{g_{0}+n-1}}&\ldots
&\int_{b_{g}}\omega _{g_{0}+n-1}
\end{array}\right),$$
$$G=\left(\begin{array}{ccc} \int_{a_{1}}\omega _{g_{0}+n}&\ldots
&\int_{a_{g_{0}}}\omega _{g_{0}+n}\\
\vdots & &\vdots \\
\int_{a_{1}}\omega _{g}&\ldots &\int_{a_{g_{0}}}\omega _{g}
\end{array}\right),$$
$$H=\left(\begin{array}{ccc} \int_{a_{g_{0}+1}}\omega
_{g_{0}+n}&\ldots
&\int_{a_{g_{0}+n-1}}\omega _{g_{0}+n}\\
\vdots & &\vdots \\
\int_{a_{g_{0}+1}}\omega _{_{g}}&\ldots
&\int_{a_{g_{0}+n-1}}\omega _{g}
\end{array}\right),$$
$$I=\left(\begin{array}{ccc} \int_{a_{g_{0}+n}}\omega
_{g_{0}+n}&\ldots
&\int_{a_{g}}\omega _{g_{0}+n}\\
\vdots & &\vdots \\
\int_{a_{g_{0}+n}}\omega _{g}&\ldots &\int_{a_{g}}\omega _{g}
\end{array}\right),$$
$$J=\left(\begin{array}{ccc} \int_{b_{1}}\omega _{g_{0}+n}&\ldots
&\int_{b_{g_{0}}}\omega _{g_{0}+n}\\
\vdots & &\vdots \\
\int_{b_{1}}\omega _{g}&\ldots &\int_{b_{g_{0}}}\omega _{g}
\end{array}\right),$$
$$K=\left(\begin{array}{ccc} \int_{b_{g_{0}+1}}\omega
_{g_{0}+n}&\ldots
&\int_{b_{g_{0}+n-1}}\omega _{g_{0}+n}\\
\vdots & &\vdots \\
\int_{b_{g_{0}+1}}\omega _{g}&\ldots &\int_{b_{g_{0}+n-1}}\omega
_{g}
\end{array}\right),$$
and $$L=\left(\begin{array}{ccc} \int_{b_{g_{0}+n}}\omega
_{g_{0}+n}&\ldots
&\int_{b_{g}}\omega _{g_{0}+n}\\
\vdots & &\vdots \\
\int_{b_{g_{0}+n}}\omega _{_{g}}&\ldots &\int_{b_{g}}\omega _{g}
\end{array}\right).$$
Notice that
\begin{eqnarray}
\int_{a_{g_{0}+1}}\omega_{j}&=&-\int_{\sigma(a_{g_{0}+1})}\omega
_{j},\nonumber\\
&=&-\int_{a_{g_{0}+1}}\sigma^{*}(\omega_{j}),\nonumber\\
&=& \left\{\begin{array}{rl}
\int_{a_{g_{0}+1}}\omega _{j},&1\leq j\leq g_{0}+n-1,\\
-\int_{a_{g_{0}+1}}\omega _{j},&g_{0}+n\leq j\leq g,
\end{array}\right.\nonumber
\end{eqnarray}
$$\vdots $$
$$
\int_{a_{g_{0}+n-1}}\omega _{j}=\left\{\begin{array}{rl}
\int_{a_{_{g_{0}+n-1}}}\omega _{j},&1\leq j\leq g_{0}+n-1,\\
-\int_{a_{_{g_{0}+n-1}}}\omega _{j},&g_{0}+n\leq j\leq g,
\end{array}\right.
$$
\begin{eqnarray}
\int_{b_{g_{0}+1}}\omega _{j}&=&-\int_{\sigma \left(
b_{g_{0}+1}\right) }\omega_{j},\nonumber\\
&=&-\int_{b_{g_{0}+1}}\sigma ^{*}\left( \omega_{j}\right), \nonumber\\
&=& \left\{\begin{array}{rl}
\int_{b_{g_{0}+1}}\omega _{j},&1\leq j\leq g_{0}+n-1,\\
-\int_{b_{g_{0}+1}}\omega _{j},&g_{0}+n\leq j\leq g,
\end{array}\right.\nonumber
\end{eqnarray}
$$\vdots $$
$$
\int_{b_{g_{0}+1}}\omega _{j}=\left\{\begin{array}{rl}
\int_{b_{g_{0}+1}}\omega _{j},&1\leq j\leq g_{0}+n-1,\\
-\int_{b_{g_{0}+1}}\omega _{j},&g_{0}+n\leq j\leq g,
\end{array}\right.
$$
\begin{eqnarray}
\int_{a_{g_{0}+n}}\omega _{j}&=&\int_{\sigma \left( a_{1}\right)
}\omega _{j},\nonumber\\
&=&\int_{a_{1}}\sigma ^{*}\left( \omega_{j}\right), \nonumber\\
&=&\left\{\begin{array}{rl}
-\int_{a_{1}}\omega _{j},&1\leq j\leq g_{0}+n-1,\\
\int_{a_{1}}\omega _{j},&g_{0}+n\leq j\leq g,
\end{array}\right.\nonumber
\end{eqnarray}
$$\vdots $$
$$
\int_{a_{g_{0}+n}}\omega _{j}=\left\{\begin{array}{rl}
-\int_{a_{_{g_{0}}}}\omega _{j},&1\leq j\leq g_{0}+n-1,\\
\int_{a_{_{g_{0}}}}\omega _{j},&g_{0}+n\leq j\leq g,
\end{array}\right.
$$ and
\begin{eqnarray}
\int_{b_{g_{0}+n}}\omega _{j}&=&\int_{\sigma \left( b_{1}\right)
}\omega _{j},\nonumber\\
&=&\int_{b_{1}}\sigma ^{*}\left( \omega_{j}\right) ,\nonumber\\
&=& \left\{\begin{array}{rl}
-\int_{b_{1}}\omega _{j},\text{ }&1\leq j\leq g_{0}+n-1,\\
\int_{b_{1}}\omega _{j},&g_{0}+n\leq j\leq g,
\end{array}\right.\nonumber
\end{eqnarray}
$$\vdots $$
$$
\int_{b_{g}}\omega _{j}=\left\{\begin{array}{rl}
-\int_{b_{g_{0}}}\omega _{j},\text{ }&1\leq j\leq g_{0}+n-1,\\
\int_{b_{g_{0}}}\omega _{j},&g_{0}+n\leq j\leq g.
\end{array}\right.
$$
Then, $$C=-A,\text{ }F=-D,\text{ }H=O,\text{ }I=G,\text{
}K=O,\text{ }L=J,$$ where $O$ is the null matrix and therefore
\begin{eqnarray}
\Omega&=&\left(\begin{array}{cccccc}
A&B&-A&D&E&-D\\
G&O&G&J&O&J
\end{array}\right),\nonumber\\
&\equiv&\left(\begin{array}{cccccc} C_1&C_2&C_3&C_4&C_5&C_6
\end{array}\right).\nonumber
\end{eqnarray}
By elementary column operation, we obtain the following matrices :
\begin{eqnarray}
\Omega_1&=&\left(\begin{array}{cccccc}
C_1&C_2&C_1+C_3&C_4&C_5&C_4+C_6
\end{array}\right),\nonumber\\
&=&\left(\begin{array}{cccccc}
A&B&O&D&E&O\\
G&O&2G&J&O&2J
\end{array}\right),\nonumber
\end{eqnarray}
\begin{eqnarray}
\nonumber\\
\Omega_2&=&\left(\begin{array}{cccccc}
C_1&C_2&C_1-C_3&C_4&C_5&C_4-C_6
\end{array}\right),\nonumber\\
&=&\left(\begin{array}{cccccc}
A&B&2A&D&E&2D\\
G&O&O&J&O&O
\end{array}\right),\nonumber
\end{eqnarray}
\begin{eqnarray}
\Omega_3&=&\left(\begin{array}{cccccc}
C_1-C_3&C_2&C_4-C_6&C_5&C_1+C_3&C_4+C_6
\end{array}\right),\nonumber\\
&=&\left(\begin{array}{cccccc}
2A&B&2D&E&O&O\\
O&O&O&O&2G&2J
\end{array}\right),\nonumber\\
&=&\left(\begin{array}{cc}
\Gamma&O\\
O&2\Delta
\end{array}\right),\nonumber
\end{eqnarray} where
\begin{eqnarray}
\Delta& =&\left(\begin{array}{cc} G&J
\end{array}\right),\\
&=&\left(\begin{array}{cccccc} \int_{a_{1}}\omega
_{g_{0}+n}&\ldots &\int_{a_{g_{0}}}
\omega _{g_{0}+n}&\int_{b_{1}}\omega _{g_{0}+n}&\ldots &\int_{b_{g_{0}}}\omega _{g_{0}+n}\\
\vdots &&\vdots &\vdots &&\vdots \\
\int_{a_{1}}\omega _{g}&\ldots &\int_{a_{g_{0}}}\omega
_{g}&\int_{b_{1}}\omega _{g}&\ldots &\int_{b_{g_{0}}}\omega _{g}
\end{array}\right),\nonumber
\end{eqnarray}
and
\begin{eqnarray}
\Gamma &=&\left(\begin{array}{cccc} 2A&B&2D&E
\end{array}\right),\\
&=&\left(\begin{array}{ccc}
2\int_{a_{1}}\omega _{1}&\ldots &2\int_{a_{1}}\omega _{g_{0}+n-1}\\
\vdots&&\vdots\\
2\int_{a_{g_{0}}}\omega _{1}&\ldots &2\int_{a_{g_{0}}}\omega _{g_{0}+n-1}\\
\int_{a_{g_{0}+1}}\omega _{1}&\ldots &\int_{a_{g_{0}+1}}\omega _{g_{0}+n-1}\\
\vdots&&\vdots\\
\int_{a_{g_{0}+n-1}}\omega _{1}&\ldots &\int_{a_{g_{0}+n-1}}\omega _{g_{0}+n-1}\\
2\int_{b_{1}}\omega _{1}&\ldots &2\int_{b_{1}}\omega _{g_{0}+n-1}\\
\vdots&&\vdots\\
2\int_{b_{g_{0}}}\omega _{1}&\ldots &2\int_{b_{g_{0}}}\omega _{g_{0}+n-1}\\
\int_{b_{g_{0}+1}}\omega _{1}&\ldots &\int_{b_{g_{0}+1}}\omega _{g_{0}+n-1}\\
\vdots&&\vdots\\
\int_{b_{g_{0}+n-1}}\omega _{1}&\ldots &\int_{b_{g_{0}+n-1}}\omega
_{g_{0}+n-1}
\end{array}\right)^\intercal.\nonumber
\end{eqnarray}
Let $$L_{\Omega
}=\{\sum_{i=1}^{g}m_{i}\int_{a_{i}}\left(\begin{array}{c}
\omega_1\\\vdots \\\omega_g
\end{array}\right )
+n_{i}\int_{b_{i}}\left(\begin{array}{c} \omega_1\\\vdots
\\\omega_g
\end{array}\right) :m_{i},n_{i}\in \mathbb{Z}\} ,$$
be the period lattice associated to $\Omega .$ Let us denote also
by $L_{\Omega _{1}},$ $L_{\Omega _{2}},$ $L_{\Omega _{3}}$ and
$L_{\Delta }$ the period lattices associated respectively to
$\Omega _{1},$ $\Omega _{2},$ $\Omega _{3}$ and $\Delta.$ Since
$$L_{\Omega }=L_{\Omega _{1}}=L_{\Omega _{2}}=L_{\Omega _{3}},$$ it
follows that the maps
$$\mathbb{C}^{g_{0}+n-1}/L_{\Gamma }:
\left(\begin{array}{c} t_1\\\vdots \\t_{g_{0}+n-1}
\end{array}\right ) \mbox{mod}L_{\Gamma }\hookrightarrow
\mathbb{C}^{g}/L_{\Omega }:\left(\begin{array}{c} t_1\\\vdots
\\t_{g_{0}+n-1}\\0\\\vdots
\\0\end{array}\right ) \mbox{mod}L_{\Omega },$$
$$\mathbb{C}^{g_{0}}/2L_{\Delta }:
\left(\begin{array}{c} t_1\\\vdots \\t_{g_{0}}
\end{array}\right ) \mbox{mod}2L_{\Delta }\hookrightarrow
\mathbb{C}^{g}/L_{\Omega }:\left(\begin{array}{c} 0\\\vdots
\\0\\t_1\\\vdots
\\t_{g_{0}}\end{array}\right ) \mbox{mod}L_{\Omega },$$
are injectives. Therefore, the tori
$\mathbb{C}^{g_{0}+n-1}/L_{\Gamma }$ and
$\mathbb{C}^{g_{0}}/2L_{\Delta }$ can be embedded into
$\mathbb{C}^{g}/L_{\Omega }$ and the map
\begin{eqnarray}
\mathbb{C}^{g}/L_{\Omega _{3}} &=&\mathbb{C}^{g_{0}+n-1}/L_{\Gamma
}\oplus \mathbb{C}^{g_{0}}/2L_{\Delta }\longrightarrow
\mathbb{C}^{g}/L_{\Omega },\nonumber\\
&&\left(\begin{array}{c} t_1\\\vdots\\t_{g}
\end{array}\right ) \mbox{mod}L_{\Omega_3}\mapsto
\left(\begin{array}{c} t_{1}\\\vdots
\\t_{g}\end{array}\right ) \mbox{mod}L_{\Omega },\nonumber
\end{eqnarray}
shows that the jacobian variety $Jac\left( \mathbb{C}_{0}\right) $
intersects the torus $\mathbb{C}^{g_{0}+n-1}/L_{\Gamma }$ in
$2^{2g_{0}}$ points. We have then the following diagramm
$$
\begin{array}{ccccccc}
\mathcal{C}&&KerN_{\varphi}&&&&\\
\null\hskip0.2cm\downarrow \sigma &&\downarrow &&&&\\
\mathcal{C}&\overset{u}{\longrightarrow}&Jac(\mathcal{C})&
\overset{u^{*}}{\longleftarrow}&Jac^*(\mathcal{C})&&\\
\null\hskip0.2cm\downarrow \varphi &&\null\hskip0.6cm\downarrow N_{\varphi}&&
\null\hskip0.4cm\uparrow N_{\varphi }^{*}&&\\
\mathcal{C}_{0}&\overset{u_{0}}{\longrightarrow}&
Jac(\mathcal{C}_{0})&\overset{u_{0}^{*}}{\longleftarrow}&Jac^*(\mathcal{C}_{0})&&\\
&&\null\hskip0.6cm\downarrow \varphi ^{*}&&&\\
&&Jac(\mathcal{C})&\subset &H^{1}(\mathcal{C},\mathcal{O}_{\mathcal{C}}^{*})&
\simeq &H^{1}(\mathcal{C}_{0},(\varphi _{*}\mathcal{O}_{\mathcal{C}})^{*})\\
&&\null\hskip0.6cm\downarrow \sigma &&&&\null\hskip0.4cm\downarrow N_{\varphi}\\
&&Jac(\mathcal{C})&\overset{N_{\varphi }}{\longrightarrow
}&Jac(\mathcal{C}_{0})&\hookrightarrow &H\text{ }^{1}(
\mathcal{C}_{0},\mathcal{O}_{\mathcal{C}_{0}}^{*})
\end{array}
$$
where
$$u:z\longmapsto \text{divisor class}(z-p) ,\text{ }p\in
\mathcal{C},\text{fixed},$$ $$u_{0}:z_{0}\longmapsto \text{divisor
class}( z_{0}-p_{0}),\text{ }p_{0}\in \mathcal{C}_{0},\text{ fixed
with}p_{0}=\varphi (p),$$
$$N_{\varphi}:Jac(\mathcal{C})\longrightarrow Jac(
\mathcal{C}_{0}),\text{ }\sum m_{i}p_{i}\longmapsto \sum
m_{i}\varphi(p_{i}),$$ is the norm mapping and
$$Jac^*(\mathcal{C})=\text{dual of } Jac(\mathcal{C}).$$ The norm
map $N_{\varphi }$ is surjective. Moreover, the dual $Jac^*\left(
\mathcal{C}\right) $ of $Jac(\mathcal{C})$ is isomorphic to
$Jac(\mathcal{C}).$ The Prym variety denoted
$Prym(\mathcal{C}/\mathcal{C}_{0})$ is defined by
$$Prym(\mathcal{C}/\mathcal{C}_{0})=\left(H^{0}\left(\mathcal{C},
\Omega_{\mathcal{C}}^{1}\right)^{-}\right)^{*}/H_{1}\left(\mathcal{C},
\mathbb{Z}\right)^{-},$$
where $-$denote the -1 eigenspace for a vector space on which j
acts. Let $D\in Prym(\mathcal{C}/\mathcal{C}_{0})$, i.e.,
\begin{eqnarray}
\left(\begin{array}{c}
0\\
\vdots \\
0\\
2\int_{0}^{D}\omega _{g_{0}+n}\\
\vdots \\
2\int_{0}^{D}\omega _{g}
\end{array}\right)\in L_{\Omega }
&\Longleftrightarrow& \left(\begin{array}{c}
\int_{0}^{D}\omega _{g_{0}+n}\\
\vdots \\
\int_{0}^{D}\omega _{g}
\end{array}\right)\in L_{\Delta },\nonumber\\
&\Longleftrightarrow& \left(\begin{array}{c}
\int_{0}^{D}\omega _{1}\\
\vdots \\
\int_{0}^{D}\omega _{g_{0}+n-1}\\
\int_{0}^{D}\omega _{g_{0}+n}\\
\vdots \\
\int_{0}^{D}\omega _{g}
\end{array}\right)=\left(\begin{array}{c}
t_{1}\\
\vdots \\
t_{g_{0}+n-1}\\
0\\
\vdots \\
0
\end{array}\right).\nonumber
\end{eqnarray}
Consequently, the Abel-Jacobi map $$Jac\left( \mathcal{C}\right)
\longrightarrow \mathcal{C}^{g}/L_{\Omega },\text{ }D\longmapsto
\left(\begin{array}{c}
\int_{0}^{D}\omega _{1}\\
\vdots \\
\int_{0}^{D}\omega _{g_{0}+n-1}\\
\int_{0}^{D}\omega _{g_{0}+n}\\
\vdots \\
\int_{0}^{D}\omega _{g}
\end{array}\right),$$
maps $Prym(\mathcal{C}/\mathcal{C}_{0}))$ bi-analtically onto
$\mathcal{C}^{g_{0}+n-1}/L_{\Gamma }.$ By Chow's theorem [9],
$Prym(\mathcal{C}/\mathcal{C}_{0})$ is thus a subabelian variety
of $Jac\left( \mathcal{C}\right) .$ More precisely, we have
\begin{eqnarray}
Prym(\mathcal{C}/\mathcal{C}_{0})&=&(KerN_{\varphi})
^{0},\nonumber\\
&=&Ker(1_{Jac(\mathcal{C})}+\sigma)
^{0},\nonumber\\
&=&Im(1_{Jac(\mathcal{C})}-\sigma )\subset Jac(
\mathcal{C}).\nonumber
\end{eqnarray}
Equivalently, the involution $\sigma $ induces an involution
$$\sigma :Jac(\mathcal{C})\longrightarrow Jac(\mathcal{C}),\text{
class of }D\longmapsto \text{class of }\sigma D,$$ and up to some
points of order two, $Jac(\mathcal{C})$ splits into an even part
$Jac\left( \mathcal{C}_{0}\right) $ and an odd part
$Prym(\mathcal{C}/\mathcal{C}_{0})$ :
$$Jac(\mathcal{C})=Prym(\mathcal{C}/\mathcal{C}_{0})\oplus
Jac\left( C_{0}\right),$$ with
\begin{eqnarray}
\dim Jac(C_{0}) &=&g_{0},\nonumber\\
\dim Jac()&=&g=2g_{0}+n-1, \nonumber\\
\dim Prym(\mathcal{C}/\mathcal{C}_{0})&=&g-g_{0}
=g_{0}+n-1.\nonumber
\end{eqnarray}
Observe that $\Delta \left(2\right)$ (resp. $\Gamma \left(
3\right)$) is the period matrix of $Jac\left(
\mathcal{C}_{0}\right)$ (resp.
$Prym(\mathcal{C}/\mathcal{C}_{0})$). Write $$\Gamma =\left(
U,\text{ }V\right) ,$$ with
$$U=\left(\begin{array}{cccccc}
2\int_{a_{1}}\omega _{1}&\ldots &2\int_{a_{g_{0}}}\omega
_{1}&\int_{a_{g_{0}+1}}\omega _{1}&\ldots
&\int_{a_{g_{0}+n-1}}\omega _{1}\\
\vdots &&\vdots &\vdots &&\vdots \\
2\int_{a_{1}}\omega _{g_{0}+n-1}&\ldots &2\int_{a_{g_{0}}}\omega
_{g_{0}+n-1}&\int_{a_{g_{0}+1}}\omega _{g_{0}+n-1}&\ldots
&\int_{a_{g_{0}+n-1}}\omega _{g_{0}+n-1}
\end{array}\right),$$
$$V=\left(\begin{array}{cccccc}
2\int_{b_{1}}\omega _{1}&\ldots &2\int_{b_{g_{0}}}\omega
_{1}&\int_{b_{g_{0}+1}}\omega _{1}&\ldots
&\int_{b_{g_{0}+n-1}}\omega _{1}\\
\vdots &&\vdots &\vdots &&\vdots\\
2\int_{b_{1}}\omega _{g_{0}+n-1}&\ldots &2\int_{b_{g_{0}}}\omega
_{g_{0}+n-1}&\int_{b_{g_{0}+1}}\omega _{g_{0}+n-1}&\ldots
&\int_{b_{g_{0}+n-1}}\omega _{g_{0}+n-1}
\end{array}\right),$$ and let us call
$$e_1=\left(\begin{array}{c}
2\int_{a_{1}}\omega _{1}\\
\vdots \\
2\int_{a_{1}}\omega _{g_{0}+n-1}
\end{array}\right),\ldots,e_{g_0}=\left(\begin{array}{cccccc}
2\int_{a_{g_{0}}}\omega _{1}\\
\vdots  \\
2\int_{a_{g_{0}}}\omega _{g_{0}+n-1}
\end{array}\right),$$
$$e_{g_{0}+1}=\left(\begin{array}{cccccc}
\int_{a_{g_{0}+1}}\omega _{1}\\
\vdots \\
\int_{a_{g_{0}+1}}
\end{array}\right),\ldots,e_{g_{0}+n-1}=\left(\begin{array}{cccccc}
\int_{a_{g_{0}+n-1}}\omega _{1}\\
\vdots \\
\int_{a_{g_{0}+n-1}}\omega _{g_{0}+n-1}
\end{array}\right).$$
Then, in the new basis $\left( \lambda _{1},\ldots ,\lambda
_{g_{0}+n-1}\right) $ where
$$\lambda _{j}=\frac{e_{j}}{\delta _{j}},\qquad\delta
_{j}=\left\{\begin{array}{rl} 1&\text{ pour }1\leq j\leq g_{0},\\
2&\text{ pour }g_{0}+1\leq j\leq g_{0}+n-1,
\end{array}\right.$$
the period matrix $\Gamma $ takes the canonical form $$(\Delta
_{\delta },\text{ }Z),$$ with $$\Delta _{\delta }=diag(\delta
_{1},\ldots ,\delta _{n})$$ and $$Z=\Delta _{\delta }U^{-1}V,$$
symmetric and $\mbox{Im}Z >0.$ Then
\begin{eqnarray}
\Gamma ^{*}&=&\left( \delta _{g_{0}+n-1}\Delta _{\delta
}^{-1},\delta _{g_{0}+n-1}\Delta _{\delta }^{-1}Z\Delta _{\delta
}^{-1}\right) ,\nonumber\\
&=&(\delta _{g_{0}+n-1}\Delta _{\delta }^{-1},\delta _{g_{0}+n-1}U^{-1}V\Delta _{\delta }^{-1}),\nonumber\\
&=&\left( \delta _{g_{0}+n-1}\Delta _{\delta }^{-1},\delta
_{g_{0}+n-1}\Delta _{\delta }^{-1}\left( U^{*}\right)
^{-1}V^{*}\right),\nonumber
\end{eqnarray} and
\begin{eqnarray}
\Gamma ^{*}&=&( U^{*},V^{*}),\nonumber\\
&=&(A\quad B\quad D\quad E),\nonumber
\end{eqnarray}
is the period matrix of the dual abelian variety
$Prym^*(\mathcal{C}/\mathcal{C}_{0}).$ The above discussion is
summed up in the following statement :

\begin{Theo}
Let $$\varphi :\mathcal{C}\longrightarrow \mathcal{C}_{0},$$ be a
double covering where $\mathcal{C}$ and $\mathcal{C}_0$ are
nonsingular algebraic curves with jacobians $Jac (\mathcal{C})$
and $Jac(\mathcal{C}_0)$. Let $$\sigma :
\mathcal{C}\longrightarrow \mathcal{C},$$ be the involution
exchanging sheets of $\mathcal{C}$ over
$\mathcal{C}_0=\mathcal{C}/\sigma.$ This involution extends by
linearity to a map (which will again be denoted by $\sigma$)
$$\sigma : Jac(\mathcal{C})\longrightarrow Jac(\mathcal{C}),$$ and
up some points of order two, $Jac(\mathcal{C})$ splits into an
even part ,i.e., $Jac (\mathcal{C})$ and an odd part (called Prym
variety) denoted $Prym(\mathcal{C}/\mathcal{C}_{0})$ and defined
by
$$Prym(\mathcal{C}/\mathcal{C}_{0})=\left( H^{0}\left(
\mathcal{C},\Omega _{\mathcal{C}}^{1}\right) ^{-}\right)
^{*}/H_{1}\left( \mathcal{C},\mathbb{Z}\right) ^{-},$$ where
$\Omega _{\mathcal{C}}^{1}$ is the sheaf of holomorphic 1-forms on
$\mathcal{C}$ and $-$ denote the -1 eigenspace for a vector space
on which j acts. To be precise, we have
$$Jac(\mathcal{C})=Prym(\mathcal{C}/\mathcal{C}_{0})\oplus
Jac\left( \mathcal{C}_{0}\right) ,$$ with $$ \dim Jac\left(
\mathcal{C}_{0}\right) =\mbox{ genus } g_{0}\mbox{ of
}\mathcal{C}_{0},$$ $$\dim Jac\left( \mathcal{C}\right) =\mbox{
genus } g\mbox{ of }\mathcal{C}=2g_{0}+n-1,$$ and $$ \dim
Prym(\mathcal{C}/\mathcal{C}_{0})=g-g_{0}=g_{0}+n+1,$$ with $2n$
branch points. Moreover, if
$$\Omega=\left(\begin{array}{cccccc}
A&B&C&D&E&F\\
G&H&I&J&K&L
\end{array}
\right),$$ is the period matrix of $Jac(\mathcal{C})$ where
$A,\ldots ,L$ denote the matrices $\left( 1\right) ,$ then the
period matrices of $Jac\left( \mathcal{C}_{0}\right) ,$
$Prym\left( \mathcal{C}/\mathcal{C}_{0}\right) $ and the dual
$Prym^*\left( \mathcal{C}/\mathcal{C}_{0}\right) $ of $Prym\left(
\mathcal{C}/\mathcal{C}_{0}\right) ,$ are respectively $$\Delta
=(G\quad H),$$ $$\Gamma =(2A \quad B \quad 2D \quad E),$$ and
$$\Gamma^*=(A \quad B \quad D \quad E).$$
\end{Theo}

\section{Algebraic complete integrability}

\subsection{A survey on abelian varieties and algebraic integrability}

We give some basic facts about integrable hamiltonian systems and
some results about abelian surfaces which will be used, as well as
the basic techniques to study two-dimensional algebraic completely
integrable systems. Let $T=\mathbb{C}/\Lambda$ be a
$n-$dimensional abelian variety where $\Lambda$ is the lattice
generated by the $2n$ columns $\lambda_1,\ldots,\lambda_{2n}$ of
the $n\times 2n$ period matrix $\Omega$ and let $D$ be a divisor
on $T.$ Define $$\mathcal{L}(\mathcal{D})=\{f \mbox{ meromorphic
on } T  : (f)\geq-\mathcal{D}\},$$ i.e., for $\mathcal{D}=\sum
k_j\mathcal{D}_j$ a function $f\in\mathcal{L}(\mathcal{D})$ has at
worst a $k_j-$fold pole along $\mathcal{D}_j.$ The divisor
$\mathcal{D}$ is called ample when a basis $(f_0,\ldots,f_N)$ of
$\mathcal{L}(k\mathcal{D})$ embeds $T$ smoothly into
$\mathbb{P}^N$ for some $k,$ via the map $$T\longrightarrow
\mathbb{P}^N,\text{ }p\longmapsto[1:f_{1}(p):...:f_{N}(p)],$$ then
$k\mathcal{D}$ is called very ample. It is known that every
positive divisor $\mathcal{D}$ on an irreducible abelian variety
is ample and thus some multiple of $\mathcal{D}$ embeds $M$ into
$\mathbb{P}^N.$ By a theorem of Lefschetz, any $k\geq 3$ will
work. Moreover, there exists a complex basis of $\mathbb{C}^n$
such that the lattice expressed in that basis is generated by the
columns of the $n\times 2n$ period matrix
$$ \left(\begin{array}{ccccc}
\delta_1&&0&|&\\
&\ddots&&|&Z\\
0&&\delta_n&|&
\end {array}\right),$$
with $Z^\top =Z, \mbox{Im}Z>0, \delta_j\in \mathbb{N}^*$ and
$\delta_j|\delta_{j+1}.$ The integers $\delta_j$ which provide the
so-called polarization of the abelian variety $M$ are then related
to the divisor as follows : \begin{equation}\label{eqn:euler} \dim
\mathcal{L}(\mathcal{D})=\delta_1\ldots \delta_n.
\end{equation}
In the case of a $2-$dimensional abelian varieties (surfaces),
even more can be stated : the geometric genus $g$ of a positive
divisor $\mathcal{D}$ (containing possibly one or several curves)
on a surface $T$ is given by the adjunction formula
\begin{equation}\label{eqn:euler}
g(\mathcal{D})=\frac{K_{T}.\mathcal{D}+\mathcal{D}.\mathcal{D}}{2}+1,
\end{equation}
where $K_T$ is the canonical divisor on $T,$ i.e., the zero-locus
of a holomorphic $2-$form, $\mathcal{D}.\mathcal{D}$ denote the
number of intersection points of $\mathcal{D}$ with
$a+\mathcal{D}$ (where $a+\mathcal{D}$ is a small translation by
$a$ of $\mathcal{D}$ on $T$), where as the Riemann-Roch theorem
for line bundles on a surface tells you that
\begin{equation}\label{eqn:euler}
\chi(\mathcal{D})=p_{a}(T)+1+\frac{1}{2}(\mathcal{D}.\mathcal{D}-\mathcal{D}K_M),
\end{equation}
where $p_{a}\left( T\right) $ is the arithmetic genus of $T$ and
$\chi(\mathcal{D})$ the Euler characteristic of $\mathcal{D}.$ To
study abelian surfaces using Riemann surfaces on these surfaces,
we recall that
\begin{eqnarray}
\chi( \mathcal{D}) &=&\dim {H}^{0}(
T,\mathcal{O}_{T}(\mathcal{D}))-\dim {H}^{1}(T,\mathcal{O}_{T}(D)),\nonumber\\
&=&\dim{\mathcal{L}}(\mathcal{D})-\dim {H}^{1}(T,\Omega ^{2}(
\mathcal{D}\otimes K_T^*)),\mbox{(Kodaira-Serre duality)},\nonumber\\
&=&\dim {\mathcal{L}}(\mathcal{D}), \mbox{(Kodaira vanishing
theorem)},
\end{eqnarray}
whenever $\mathcal{D}\otimes K_T^*$ defines a positive line
bundle. However for abelian surfaces, $K_T$ is trivial and
$p_a(T)=-1;$ therefore combining relations (4), (5), (6) and (7),
we obtain
\begin{eqnarray}
\chi \left( \mathcal{D}\right)&=&\dim
{\mathcal{L}}(\mathcal{D}),\nonumber\\
&=&\frac{\mathcal{D}.\mathcal{D}}{2},\nonumber\\
&=&g\left(\mathcal{D}\right)-1,\nonumber\\
&=&\delta_1 \delta_2.
\end{eqnarray}

Let $M$ be a $2n$-dimensional differentiable manifold and $\omega
$ a closed non-degenerate differential $2$-form. The pair $\left(
M,\omega \right) $ is called a symplectic manifold. Let
$H:M\longrightarrow \mathbb{R}$ be a smooth function. A
hamiltonian system on $\left( M,\omega \right) $ with hamiltonian
$H$ can be written in the form
\begin{equation}\label{eqn:euler}
\dot q_{1}=\frac{\partial H}{\partial p_{1}},...,\dot
q_{n}=\frac{\partial H}{\partial p_{n}},\quad \dot
p_{1}=-\frac{\partial H}{\partial q_{1}},...,\dot
p_{n}=-\frac{\partial H}{\partial q_{n}},
\end{equation}
where $\left(q_{1},...,q_{n},p_{1},...,p_{n}\right) $ are
coordinates in $M$. Thus the hamiltonian vector field $X_{H}$ is
defined by $$X_{H}=\sum_{k=1}^{n}\left( \frac{\partial H}{\partial
p_{k}}\frac{\partial }{\partial q_{k}}-\frac{\partial H}{\partial
q_{k}}\frac{\partial }{\partial p_{k}}\right).$$ If $F$ is a
smooth function on the manifold $M,$ the Poisson bracket $\left\{
F,H\right\} $ of $F$ and $H$ is defined by
\begin{eqnarray}
X_{H}F&=&\sum_{k=1}^{n}\left( \frac{\partial H}{\partial
p_{k}}\frac{\partial F}{\partial q_{k}}-\frac{\partial H}{\partial
q_{k}}\frac{\partial F}{\partial p_{k}}\right),\nonumber\\
&=&\left\{ F,H\right\}.
\end{eqnarray}
A function $F$ is an invariant (first integral) of the hamiltonian
system (9) if and only if the Lie derivative of $F$ with respect
$X_{H}$ is identically zero. The functions $F$ and $H$ are said to
be in involution or to commute, if $\left\{F,H\right\} =0.$ Note
that equations (9) and (10) can be written in more compact form
$$
\dot{x}=J\frac{\partial H}{\partial x}\text{ },\text{\quad
}x=\left( q_{1},...,q_{n},p_{1},...,p_{n}\right) ^{\top },
$$
\begin{eqnarray}
\left\{ F,H\right\}&=&\left\langle \frac{\partial F}{\partial x},
J\frac{\partial H}{\partial x}\right\rangle,\nonumber\\
&=&\sum_{k,l=1}^{n}J_{kl}\frac{\partial F}{\partial
x_{k}}\frac{\partial H}{\partial x_{l}},\nonumber
\end{eqnarray}
with
$J=\left[\begin{array}{cc}
O&I\\
-I&O
\end{array}\right],$
a skew-symmetric matrix where $I$ is the $n\times n$ unit matrix
and $O$ the $n\times n$ zero matrix. A hamiltonian system is
completely integrable in the sense of Liouville if there exist $n$
invariants $H_{1}=H,H_{2},\ldots ,H_{n}$ in involution (i.e., such
that the associated Poisson bracket $\{H_{k},H_{l}\}$ all vanish)
with linearly independent gradients (i.e., $dH_{1}\wedge ...\wedge
dH_{n}\neq 0$). For generic $c=\left( c_{1},...,c_{n}\right) $ the
level set
$$
M_c=\{H_{1}=c_{1},\ldots ,H_{n}=c_{n}\},
$$
will be an $n$-manifold, and since
$$X_{H_{k}}H_{l}=\{H_{k},H_{l}\}=0,$$ the integral curves of each
$X_{H_{k}}$ will lie in $M_c$ and the vector fields $X_{H_{k}}$
span the tangent space of $M_c$. By a theorem of Arnold [4,15], if
$M_c$ is compact and connected, it is diffeomorphic to an
$n$-dimensional real torus and each vector field will define a
linear flow there. To be precise, in some open neighbourhood of
the torus one can introduce regular symplectic coordinates
$s_{1},\ldots ,s_{n},\varphi _{1},\ldots ,\varphi _{n}$ in which
$\omega $ takes the canonical form $$\omega
=\sum_{k=1}^{n}ds_{k}\wedge d\varphi _{k}.$$ Here the functions
$s_{k}$ (called action-variables) give coordinates in the
direction transverse to the torus and can be expressed
functionally in terms of the first integrals $H_{k}.$ The
functions $\varphi _{k}$ (called angle-variables) give standard
angular coordinates on the torus, and every vector field
$X_{H_{k}}$ can be written in the form $$\dot{\varphi
_k}=h_{k}\left( s_{1},\ldots ,s_{n}\right),$$ that is, its
integral trajectories define a conditionally-periodic motion on
the torus. Consequently, in a neighbourhood of the torus the
hamiltonian vector field $X_{H_{k}}$ take the following form
$$\dot{s_k}=0,\quad \dot{\varphi _k}=h_{k}\left( s_{1},\ldots
,s_{n}\right),$$ and can be solved by quadratures.

Consider now hamiltonian problems of the form
\begin{equation}\label{eqn:euler}
X_{H}:\dot {x}=J\frac{\partial H}{\partial x}\text{ }\equiv
f(x),\text{ }x\in \mathbb{R}^{m},
\end{equation}
where $H$\ is the hamiltonian and $J=J(x)$\ is a skew-symmetric
matrix with polynomial entries in $x,$\ for which the
corresponding Poisson bracket $$\{H_{i},H_{j}\}=\langle
\frac{\partial H_i}{\partial x},J\frac{\partial H_j} {\partial
x}\rangle,$$ satisfies the Jacobi identities. The system (11) with
polynomial right hand side will be called algebraically completely
integrable (a.c.i.) in the sense of Adler-van Moerbeke [1] when :\\
\null\hskip0.4cm $(i)$\ The system possesses $n+k$\ independent
polynomial invariants $H_{1},...,H_{n+k}$ of which $k$ lead to
zero vector fields $$J\frac{\partial H_{n+i}}{\partial x}\left(
x\right) =0,\quad1\leq i\leq k,$$ the $n$ remaining ones are in
involution (i.e., $\left\{ H_{i},H_{j}\right\} =0$) and $m=2n+k.$\
For most values of $c_{i}\in \mathbb{R},$ the invariant varieties
$\overset{n+k}{\underset{i=1}{\bigcap }}\left\{ x\in
\mathbb{R}^{m}:H_{i}=c_{i}\right\} $\ are assumed compact and
connected. Then, according to the Arnold-Liouville theorem, there
exists a diffeomorphism
$$\overset{n+k}{\underset{i=1}{\bigcap }}\left\{ x\in
\mathbb{R}^{m}:H_{i}=c_{i}\right\} \rightarrow
\mathbb{R}^{n}/Lattice,$$ and the solutions of the system (11) are
straight lines motions on these tori.\\
\null\hskip0.4cm $(ii)$\ The invariant varieties, thought of as
affine varieties
$$M_c=\overset{n+k}{\underset{i=1}{\bigcap }}\{H_{i}=c_{i},x\in
\mathbb{C}^{m}\},$$ in $\mathbb{C}^{m}$ can be completed into
complex algebraic tori, i.e., $$ M_c\cup
\mathcal{D}=\mathbb{C}^{n}/Lattice,$$ where
$\mathbb{C}^{n}/Lattice$ is a complex algebraic torus (i.e.,
abelian variety) and $\mathcal{D}$ a divisor.

\begin{Rem}
Algebraic means that the torus\ can be defined as an intersection
$$\displaystyle{\bigcap_i\{ P_{i}(X_{0},...,X_{N})=0\}},$$
involving a large number of homogeneous polynomials $P_{i}.$\ In
the natural coordinates $(t_{1},...,t_{n})$\ of
$\mathbb{C}^{n}$/$Lattice$\ coming from $\mathbb{C}^{n},$\ the
functions $x_{i}=x_{i}(t_{1},...,t_{n})$\ are meromorphic and (11)
defines straight line motion on $\mathbb{C}^{n}/Lattice.$\
Condition $(i)$\ means, in particular, there is an algebraic map
$$(x_{1}(t),...,x_{m}(t) ) \longmapsto (\mu_{1}(t)
,...,\mu_{n}(t)),$$ making the following sums linear in $t$\
:$$\sum_{i=1}^{n}\int_{\mu_{i}(0)}^{\mu_{i}(t)}\omega _{j}=d
_{j}t\text{ },\text{ }1\leq j\leq n,\text{ }d _{j}\in
\mathbb{C},$$where $\omega _{1},...,\omega _{n}$\ denote
holomorphic differentials on some algebraic curves.
\end{Rem}

Adler and van Moerbeke $[1] $\ have shown that the existence of a
coherent set of Laurent solutions :
\begin{equation}\label{eqn:euler}
x_{i}=\sum_{j=0}^{\infty }x_{i}^{(j) }t^{j-k_{i}},\text{\quad
}k_{i}\in \mathbb{Z},\text{ \quad some }k_{i}> 0,
\end{equation}
depending on $dim\ (phase$\ $space)-1=$\ $m-1$\ free parameters is
necessary and sufficient for a hamiltonian system with the right
number of constants of motion to be a.c.i. So, if the hamiltonian
flow (11) is a.c.i., it means that the variables $x_{i}$\ are
meromorphic on the torus $\mathbb{C}^{n}/Lattice$\ and by
compactness they must blow up along a codimension one subvariety
(a divisor) $\mathcal{D}\subset \mathbb{C}^{n}/Lattice.$\ By the
a.c.i. definition, the flow (11) is a straight line motion in
$\mathbb{C}^{n}/Lattice$\ and thus it must hit the divisor
$\mathcal{D}$\ in at least one place. Moreover through every point
of $\mathcal{D},$ there is a straight line motion and therefore a
Laurent expansion around that point of intersection. Hence the
differential equations must admit Laurent expansions which depend
on the $n-1$\ parameters defining $\mathcal{D}$\ and the $n+k$\
constants $c_{i}$\ defining the torus $\mathbb{C}^{n}/Lattice$\ ,
the total count is therefore $m-1=dim\ (phase\ space)-1$\
parameters.

Assume now hamiltonian flows to be (weight)-homogeneous with a
weight $\nu _{i}\in \mathbb{N},$ going with each variable $x_{i},$
i.e.,
$$f_{i}\left( \alpha ^{\nu _{1}}x_{1},...,\alpha ^{\nu _{m}}x_{m}\right)
=\alpha ^{\nu _{i}+1}f_{i}\left( x_{1},...,x_{m}\right) , \text{
}\forall \alpha \in \mathbb{C}.$$ Observe that then the constants
of the motion $H$ can be chosen to be (weight)-homogeneous :
$$H\left( \alpha ^{\nu _{1}}x_{1},...,\alpha ^{\nu _{m}}x_{m}\right)=
\alpha ^{k}H\left( x_{1},...,x_{m}\right) ,\text{ }k\in
\mathbb{Z}.$$ If the flow is algebraically completely integrable,
the differential equations $\left(11\right) $ must admits Laurent
series solutions $\left( 12\right) $ depending on $m-1$ free
parameters. We must have $k_{i}=\nu _{i}$ and coefficients in the
series must satisfy at the 0$^{th}$step non-linear equations,

\begin{equation}\label{eqn:euler}
f_{i}\left( x_{1}^{\left( 0\right) },...,x_{m}^{\left( 0\right)
}\right) +g_{i}x_{i}^{\left( 0\right) }=0,\text{ }1\leq i\leq m,
\end{equation}
and at the k$^{th}$step, linear systems of equations :
\begin{equation}\label{eqn:euler}
\left( L-kI\right) z^{\left( k\right) }=
\left\{\begin{array}{rl} 0&\mbox{ for } k=1\\
\mbox{some polynomial in}& x^{\left( 1\right) },...,x^{\left(
k-1\right)} \mbox{ for } k>1,
\end{array}\right.
\end{equation}
where
\begin{eqnarray}
L&=&\text{ Jacobian map of }\left( 13\right),\nonumber\\
&=&\text{ }\frac{\partial f}{\partial z}+gI\mid _{z=z^{\left(
0\right)}}.\nonumber
\end{eqnarray}
If $m-1$ free parameters are to appear in the Laurent series, they
must either come from the non-linear equations $\left( 13\right) $
or from the eigenvalue problem $\left( 14\right) ,$ i.e., $L$ must
have at least $m-1$ integer eigenvalues. These are much less
conditions than expected, because of the fact that the homogeneity
$k$ of the constant $H$ must be an eigenvalue of $L$ Moreover the
formal series solutions are convergent as a consequence of the
majorant method $\left[ 1\right].$ So, the question is how does
one prove directly that the system is effectively a.c.i. with
abelian space coordinates? The idea of the direct proof used by
Adler and van Moerbeke is closely related to the geometric spirit
of the (real) Arnold-Liouville theorem discuted above. Namely, a
compact complex $n$-dimensional variety on which there exist $n$
holomorphic commuting vector fields which are independent at every
point is analytically isomorphic to a $n$-dimensional complex
torus $\mathbb{C}^{n}/Lattice$ and the complex flows generated by
the vector fields are straight lines on this complex torus. Now,
the affine variety $M_c$ is not compact and the main problem will
be to complete $M_c$\ into a non singular compact complex
algebraic variety $\widetilde{M_c}=M_c\cup \mathcal{D}$ in such a
way that the vector fields $X_{H_{1}},...,X_{H_{m}}$ extend
holomorphically along the divisor $\mathcal{D}$ and remain
independent there. If this is possible, $\widetilde{M_c}$ is an
algebraic complex torus, i.e., an abelian variety and the
coordinates $x_{i}$ restricted to $M_c$ are abelian functions. A
naive guess would be to take the natural compactification
$\overline{M_c}$ of $M_c$ in
$\mathbb{\mathbb{P}}^{m}(\mathbb{C}).$ Indeed, this can never work
for a general reason: an abelian variety $\widetilde{M_c}$ of
dimension bigger or equal than two is never a complete
intersection, that is it can never be descriped in some projective
space $\mathbb{P}^{m}( \mathbb{C}) $ by $m$-dim $\widetilde{M_c}$
global polynomial homogeneous equations. In other words, if $M_c$
is to be the affine part of an abelian variety, $\overline{M_c}$
must have a singularity somewhere along the locus at infinity
$\overline{M_c}\cap \{ X_{0}=0\}.$ In fact, Adler and van Moerbeke
[1] showed that the existence of meromorphic solutions to the
differential equations $(2) $ depending on $m-1$ free parameters
can be used to manufacture the tori, without ever going through
the delicate procedure of blowing up and down. Information about
the tori can then be gathered from the divisor. An exposition of
such methods and their applications can be found in [1,2].

\subsection{The Hénon-Heiles system}
The Hénon-Heiles system
\begin{equation}\label{eqn:euler}
\dot q_{1}=\frac{\partial H}{\partial p_{1}},\quad \dot
q_{2}=\frac{\partial H}{\partial p_{2}},\quad \dot
p_{1}=-\frac{\partial H}{\partial q_{1}},\quad \dot
p_{2}=-\frac{\partial H}{\partial q_{2}},
\end{equation}
with $$H\equiv H_{1}=\frac{1}{2}\left(
p_{1}^{2}+p_{2}^{2}+aq_{1}^{2}+bq_{2}^{2}\right)
+q_{1}^{2}q_{2}+6q_{2}^{3},$$ has another constant of motion
$$H_{2}=q_{1}^{4}+4q_{1}^{2}q_{2}^{2}-4p_{1}\left( p_{1}q_{2}-p_{2}q_{1}\right)
+4aq_{1}^{2}q_{2}+\left( 4a-b\right) \left(
p_{1}^{2}+aq_{1}^{2}\right),$$ where $a,$ $b,$ are constant
parameters and $q_{1},q_{2},p_{1},p_{2}$ are canonical coordinates
and momenta, respectively. First studied as a mathematical model
to describe the chaotic motion of a test star in an axisymmetric
galactic mean gravitational field this system is widely explored
in other branches of physics. It well-known from applications in
stellar dynamics, statistical mechanics and quantum mechanics. It
provides a model for the oscillations of atoms in a three-atomic
molecule. The system $\left(15\right) $ possesses Laurent series
solutions depending on $3$ free parameters $\alpha ,\beta ,\gamma
,$ namely
\begin{eqnarray}
q_{1}&=&\frac{\alpha }{t}+\left( \frac{\alpha
^{3}}{12}+\allowbreak \frac{\alpha A}{2}-\frac{\alpha
B}{12}\right) t+\beta t^{2}+q_{1}^{\left( 4\right)
}t^{3}+q_{1}^{\left( 5\right)
}t^{4}+q_{1}^{\left( 6\right) }t^{5}+\cdots,\nonumber\\
q_{2}&=&-\frac{1}{t^{2}}+\frac{\alpha
^{2}}{12}-\frac{B}{12}+\left( \frac{\alpha ^{4}}{48}+\frac{\alpha
^{2}A}{10}-\frac{\alpha ^{2}B}{60}-\frac{B^{2}}{240}\right)
t^{2}+\frac{\alpha \beta }{3}t^{3}+\gamma t^{4}+\cdots,\nonumber
\end{eqnarray}
where $p_{1}=\overset{.}{q}_{1},\text{ }p_{2}=\overset{.}{q}_{2}$
and
\begin{eqnarray}
q_{1}^{\left( 4\right) }&=&\frac{\alpha AB}{24}-\frac{\alpha
^{5}}{72}+\frac{11\alpha ^{3}B}{720}-\frac{11\alpha
^{3}A}{120}-\frac{\alpha B^{2}}{720}-\frac{\alpha
A^{2}}{8},\nonumber\\
q_{1}^{\left( 5\right) }&=&-\frac{\beta \alpha
^{2}}{12}+\frac{\beta B}{60}-\frac{A\beta }{10},\nonumber\\
q_{1}^{\left( 6\right) }&=&-\frac{\alpha \gamma }{9}-\frac{\alpha
^{7}}{15552}-\frac{\alpha ^{5}A}{2160}+\frac{\alpha
^{5}B}{12960}+\frac{\alpha ^{3}B^{2}}{25920}+\frac{\alpha
^{3}A^{2}}{1440}-\frac{\alpha ^{3}AB}{4320}+\frac{\alpha
AB^{2}}{1440}\nonumber\\
&&-\frac{\alpha B^{3}}{19440}-\frac{\alpha
A^{2}B}{288}+\frac{\alpha A^{3}}{144}.\nonumber
\end{eqnarray}
Let $\mathcal{D}$ be the pole solutions restricted to the surface
$$M_{c}=\overset{2}{\underset{i=1}{\bigcap }}\left\{ x\equiv
(q_1,q_2,p_1,p_2)\in \mathbb{C}^{4}, H_{i}\left( x\right)
=c_{i}\right\},$$ to be precise $\mathcal{D}$ is the closure of
the continuous components of the set of Laurent series solutions
$x\left( t\right)$ such that $H_{i}\left( x\left( t\right) \right)
=c_{i},\text{ }1\leq i\leq 2$, i.e., $
\mathcal{D}=t^{0}-\text{coefficient of }M_c. $ Thus we find an
algebraic curve defined by
\begin{equation}\label{eqn:euler}
\mathcal{D} : \beta ^{2}=P_{8}(\alpha),
\end{equation}
where $$
P_{8}\left( \alpha \right) =-\frac{7}{15552}\alpha ^{8}-\frac{1}{432}%
\left( 5A-\frac{13}{18}B\right) \alpha ^{6}
-\frac{1}{36}\left( \frac{671}{15120}B^{2}+\frac{17}{7}A^{2}-\frac{943}{%
1260}BA\right) \alpha ^{4} $$
$$-\frac{1}{36}\left( 4A^{3}-\frac{1}{2520}B^{3}-\frac{13}{6}A^{2}B+\frac{2}{%
9}AB^{2}-\frac{10}{7}c_{1}\right) \alpha ^{2}+\frac{1}{36}c_{2}.
$$ The curve $\mathcal{D}$ determined by an eight-order equation
is smooth, hyperelliptic and its genus is $3$. Moreover, the map
\begin{equation}\label{eqn:euler}
\sigma :\mathcal{D}\longrightarrow \mathcal{D},\text{ }(\beta
,\alpha ) \longmapsto (\beta ,-\alpha ),
\end{equation}
is an involution on $\mathcal{D}$ and the quotient
$\mathcal{E}=\mathcal{D}/\sigma $ is an elliptic curve defined by
\begin{equation}\label{eqn:euler}
\mathcal{E}: \beta ^{2}=P_{4}(\zeta),
\end{equation}
where $P_{4}\left( \zeta \right) $ is the degree $4$ polynomial in
$\zeta =\alpha ^{2}$ obtained from $\left( 16\right) .$ The
hyperelliptic curve $\mathcal{D}$ is thus a $2$-sheeted ramified
covering of the elliptic curve $\mathcal{E}\left( 18\right) ,$
\begin{equation}\label{eqn:euler}
\rho: \mathcal{D}\longrightarrow \mathcal{E},\text{ }(\beta
,\alpha )\longmapsto (\beta ,\zeta ),
\end{equation}
ramified at the four points covering $\zeta =0$ and $\infty .$
Following the methods in Adler-van Moerbeke [1,2], the affine
surface $M_c$ completes into an abelian surface $\widetilde{M}_c$,
by adjoining the divisor $\mathcal{D}$. The latter defines on
$\widetilde{M}_c$ a polarization $(1,2).$ The divisor
$2\mathcal{D}$ is very ample and the functions $$1,\text{
}y_{1},\text{ }y_{1}^{2},\text{ }y_{2},\text{ }x_{1}, \text{
}x_{1}^{2}+y_{1}^{2}y_{2},\text{ }x_{2}y_{1}-2x_{1}y_{2},\text{ }
x_{1}x_{2}+2Ay_{1}y_{2}+2y_{1}y_{2}^{2},$$ embed $\widetilde{M}_c$
smoothly into $\mathbb{CP}^7$ with polarization $(2,4),$ according
to (8). Then the system (15) is algebraic complete integrable and
the corresponding flow evolues on an abelian surface
$\widetilde{M}_c = \mathbb{C}^2 /\mbox{lattice},$ where the
lattice is generated by the period matrix $\left(
\begin{array}{llll}
2 & 0 & a & c \\
0 & 4 & c & b
\end{array}
\right)$, $\text{ Im}\left(
\begin{array}{ll}
a & c \\
c & b
\end{array}
\right) >0.$

\begin{Theo}
The abelian surface $\widetilde{M}_{c}$ which completes the affine
surface $M_{c}$ is the dual Prym variety $Prym^{*}\left(
\mathcal{D}/\mathcal{E}\right) $ of the genus $3$ hyperelliptic
curve $\mathcal{D}$ (16) for the involution $\sigma$ interchanging
the sheets of the double covering $\rho $ (19) and the problem
linearizes on this variety.
\end{Theo}
\emph{Proof}. Let $\left(
a_{1},a_{2},a_{3},b_{1},b_{2},b_{3}\right) $ be a canonical
homology basis of $\mathcal{D}$ such that $$\sigma \left(
a_{1}\right)=a_{3},\quad\sigma \left(
b_{1}\right)=b_{3},\quad\sigma \left(
a_{2}\right)=-a_{2},\quad\sigma \left( b_{2}\right)=-b_{2},$$ for
the involution $\sigma $ (17). As a basis of holomorphic
differentials $\omega _{0},\omega _{1},\omega _{2}$ on the curve
$\mathcal{D}$ (16) we take the differentials
$$\omega _{1}=\frac{\alpha ^{2}d\alpha }{\beta },\text{ }\omega _{2}=\frac{d\alpha }{\beta },
\text{ }\omega _{3}=\frac{\alpha d\alpha }{\beta},$$ and obviously
$$\sigma ^{*}(\omega _{1})=-\omega _{1},\quad\sigma ^{*}(\omega
_{2})=-\omega _{2},\quad\sigma ^{*}(\omega _{3})=\omega _{3}.$$ We
consider the period matrix $\Omega$ of $Jac(\mathcal{D})$
$$\Omega=\left(
\begin{array}{cccccc}
\int_{a_{1}}\omega _{1}&\int_{a_{2}}\omega _{1}&\int_{a_{3}}\omega
_{1}&\int_{b_{1}}\omega _{1}&\int_{b_{2}}\omega _{1}&\int_{b_{3}}\omega _{1}\\
\int_{a_{1}}\omega _{2}&\int_{a_{2}}\omega _{2}&\int_{a_{3}}\omega
_{2}&\int_{b_{1}}\omega _{2}&\int_{b_{2}}\omega _{2}&\int_{b_{3}}\omega _{2}\\
\int_{a_{1}}\omega _{3}&\int_{a_{2}}\omega _{3}&\int_{a_{3}}\omega
_{3}&\int_{b_{1}}\omega _{3}&\int_{b_{2}}\omega
_{3}&\int_{b_{3}}\omega _{3}
\end{array}
\right).
$$
By theorem 1,
$$\Omega=\left(
\begin{array}{cccccc}
\int_{a_{1}}\omega _{1}&\int_{a_{2}}\omega
_{1}&-\int_{a_{1}}\omega _{1}
&\int_{b_{1}}\omega _{1}&\int_{b_{2}}\omega _{1}&-\int_{b_{1}}\omega _{1}\\
\int_{a_{1}}\omega _{2}&\int_{a_{2}}\omega
_{2}&-\int_{a_{1}}\omega _{2}
&\int_{b_{1}}\omega _{2}&\int_{b_{2}}\omega _{2}&-\int_{b_{1}}\omega _{2}\\
\int_{a_{1}}\omega _{3}&0&\int_{a_{1}}\omega
_{3}&\int_{b_{1}}\omega _{3}&0&\int_{b_{1}}\omega _{3}
\end{array}
\right),
$$
and therefore the period matrices of $Jac(\mathcal{E})$(i.e.,
$\mathcal{E}$), $Prym(\mathcal{D}/\mathcal{E})$ and
$Prym^*(\mathcal{D}/\mathcal{E})$ are respectively
$$\Delta=(\int_{a_{1}}\omega _{3}\quad \int_{b_{1}}\omega _{3}),$$
$$
\Gamma=\left(
\begin{array}{cccc}
2\int_{a_{1}}\omega _{1}&\int_{a_{2}}\omega
_{1}&2\int_{b_{1}}\omega _{1}
&\int_{b_{2}}\omega _{1}\\
2\int_{a_{1}}\omega _{2}&\int_{a_{2}}\omega
_{2}&2\int_{b_{1}}\omega _{2} &\int_{b_{2}}\omega _{2}
\end{array}
\right),
$$ and
$$
\Gamma^*=\left(
\begin{array}{cccc}
\int_{a_{1}}\omega _{1}&\int_{a_{2}}\omega _{1}&\int_{b_{1}}\omega
_{1}
&\int_{b_{2}}\omega _{1}\\
\int_{a_{1}}\omega _{2}&\int_{a_{2}}\omega _{2}&\int_{b_{1}}\omega
_{2} &\int_{b_{2}}\omega _{2}
\end{array}
\right).
$$
Let $$L_{\Omega
}=\{\sum_{i=1}^{3}m_{i}\int_{a_{i}}\left(\begin{array}{c}
\omega_1\\\omega_2 \\\omega_3
\end{array}\right )
+n_{i}\int_{b_{i}}\left(\begin{array}{c} \omega_1\\\omega_2
\\\omega_3
\end{array}\right) :m_{i},n_{i}\in \mathbb{Z}\} ,$$
be the period lattice associated to $\Omega.$ Let us denote also
by $L_{\Delta},$ the period lattice associated $\Delta.$ According
to theorem 1, we get the following diagram
$$
\begin{array}{ccccccccc}
&&&&0&&&\\
&&&&\downarrow&&&&\\
&&&&\mathcal{E}&&\mathcal{D}&&\\
&&&&\quad \downarrow \varphi^*&\swarrow&\quad \downarrow \varphi&&\\
0&\longrightarrow&\ker
N_\varphi&\longrightarrow&Prym(\mathcal{D}/\mathcal{E}) \oplus
\mathcal{E}=Jac(\mathcal{D})&\overset{N_{\varphi
}}{\longrightarrow }
&\mathcal{E}&\longrightarrow&0\\
&&&\searrow \tau&\downarrow&&&&\\
&&&&\widetilde{M}_c=M_c\cup 2\mathcal{D}\simeq \mathbb{C}^2/\mbox{lattice}&&&&\\
&&&&\downarrow&&&&\\
&&&&0&&&&
\end{array}
$$
The polarization map $$\tau :
Prym(\mathcal{D}/\mathcal{E})\longrightarrow
\widetilde{M}_{c}=Prym^{*}(\mathcal{D}/\mathcal{E}),$$ has kernel
$(\varphi ^{*}\mathcal{E})\simeq \mathbb{Z}_{2}\times
\mathbb{Z}_{2}$ and the induced polarization on
$Prym(\mathcal{D}/\mathcal{E})$ is of type (1,2). Let
$$\widetilde{M}_{c}\longrightarrow \mathbb{C}^{2}/L_{\Lambda}:p\longmapsto
\int_{p_{0}}^{p}\binom{dt_{1}}{dt_{2}},$$ be the uniformizing map
where $dt_1, dt_2$ are two differentials on $\widetilde{M}_c$
corresponding to the flows generated respectively by $H_1, H_2$
such that : $$dt_1|_\mathcal{D}=\omega_1,\qquad
dt_2|_\mathcal{D}=\omega_2,$$ and
$$L_{\Lambda
}=\{\sum_{k=1}^{4}n_{k}\left(\begin{array}{c}
\int_{\nu_{k}}dt_1\\\int_{\nu_{k}}dt_2
\end{array}\right ):n_{k}\in \mathbb{Z}\} ,$$
is the lattice associated to the period matrix
$$
\Lambda=\left(
\begin{array}{cccc}
\int_{\nu_{1}}dt _{1}&\int_{\nu_{2}}dt _{1}&\int_{\nu_{4}}dt _{1}
&\int_{\nu_{4}}dt _{1}\\
\int_{\nu_{1}}dt _{2}&\int_{\nu_{2}}dt _{2}&\int_{\nu_{3}}dt _{2}
&\int_{\nu_{4}}dt _{2}
\end{array}
\right),
$$
where $(\nu_{1},\nu_{2},\nu_{3},\nu_{4})$ is a basis of
$H_{1}(\widetilde{M}_{c},\mathbb{Z})$. By the Lefschetz theorem on
hyperplane section [9], the map
$$H_{1}(\mathcal{D},\mathbb{Z})\longrightarrow
H_{1}(\widetilde{M}_{c},\mathbb{Z}),$$ induced by the inclusion
$\mathcal{D}\hookrightarrow $ $\widetilde{M}_{c}$\ is surjective
and consequently we can find $4$\ cycles $\nu _{1},\nu _{2},\nu
_{3},\nu _{4}$\ on the curve $D$ such that
$$
\Lambda=\left(
\begin{array}{cccc}
\int_{\nu_{1}}\omega _{1}&\int_{\nu_{2}}\omega
_{1}&\int_{\nu_{4}}\omega _{1}
&\int_{\nu_{4}}\omega _{1}\\
\int_{\nu_{1}}\omega _{2}&\int_{\nu_{2}}\omega
_{2}&\int_{\nu_{3}}\omega _{2} &\int_{\nu_{4}}\omega _{2}
\end{array}
\right),
$$
and $$L_{\Lambda }=\{\sum_{k=1}^{4}n_{k}\left(\begin{array}{c}
\int_{\nu_{k}}\omega_1\\\int_{\nu_{k}}\omega_2
\end{array}\right ):n_{k}\in \mathbb{Z}\}.$$
The  cycles $\nu _{1},\nu _{2},\nu _{3},\nu _{4}$ in $\mathcal{D}$
which we look for are $a_{1},b_{1},a_{2},b_{2}$ and they generate
$H_{1}(\widetilde{M}_{c},\mathbb{Z})$ such that
$$
\Lambda=\left(
\begin{array}{cccc}
\int_{a_{1}}\omega _{1}&\int_{b_{1}}\omega _{1}&\int_{a_{2}}\omega
_{1}&\int_{b_{2}}\omega _{1}\\
\int_{a_{1}}\omega _{2}&\int_{b_{1}}\omega _{2}&\int_{a_{2}}\omega
_{2} &\int_{b_{2}}\omega _{2}
\end{array}
\right),
$$
is a Riemann matrix. We show that $\Lambda=\Gamma^*$ ,i.e., the
period matrix of $Prym^*(\mathcal{D}/\mathcal{E})$ dual of
$Prym(\mathcal{D}/\mathcal{E})$. Consequently $\widetilde{M}_{c}$
and $Prym^*(\mathcal{D}/\mathcal{E})$ are two abelian varieties
analytically isomorphic to the same complex torus
$\mathbb{C}^{2}/L_{\Lambda}.$ By Chow's theorem,
$\widetilde{M}_{c}$\ and $Prym^*(\mathcal{D}/\mathcal{E})$ are
then algebraically isomorphic.

\subsection{The Kowalewski rigid body motion}

The motion for the Kowalewski's top is governed by the equations
\begin{eqnarray}
\overset{.}{m}&=&m\wedge \lambda m+\gamma \wedge l,\\
\overset{.}{\gamma }&=&\gamma \wedge \lambda m,\nonumber
\end{eqnarray}
where $m,\gamma $ and $l$ denote respectively the angular
momentum, the directional cosine of the $z$-axis (fixed in space),
the center of gravity which after some rescaling and normalization
may be taken as $l=\left( 1,0,0\right) $ and $\lambda m=\left(
m_{1}/2,m_{2}/2,m_{3}/2\right) .$ The system (20) can be written
\begin{eqnarray}
\overset{.}{m}_{1}&=&m_{2}m_{3},\nonumber\\
\overset{.}{m}_{2}&=&-m_{1}m_{3}+2\gamma
_{3},\nonumber\\
\overset{.}{m}_{3}&=&-2\gamma_{2},\\
\overset{.}{\gamma
}_{1}&=&2m_{3}\gamma _{2}-m_{2}\gamma _{3},\nonumber\\
 \overset{.}{\gamma }_{2}&=&m_{1}\gamma
_{3}-2m_{3}\gamma _{1},\nonumber\\
\overset{.}{\gamma }_{3}&=&m_{2}\gamma _{1}-m_{1}\gamma
_{2},\nonumber
\end{eqnarray}
with constants of motion
\begin{eqnarray}
H_{1}&=&\frac{1}{2}\left( m_{1}^{2}+m_{2}^{2}\right) +m_{3}^{2}+2\gamma _{1}=c_{1},\nonumber\\
H_{2}&=&m_{1}\gamma _{1}+m_{2}\gamma _{2}+m_{3}\gamma _{3}=c_{2},\\
H_{3}&=&\gamma _{1}^{2}+\gamma _{2}^{2}+\gamma _{3}^{2}=c_{3}=1,\nonumber\\
H_{4}&=&\left( \left( \frac{m_{1}+im_{2}}{2}\right) ^{2}-\left(
\gamma _{1}+i\gamma _{2}\right) \right) \left( \left(
\frac{m_{1}-im_{2}}{2}\right) ^{2}-\left( \gamma _{1}-i\gamma
_{2}\right) \right) =c_{4}.\nonumber
\end{eqnarray}
The system (21) admits two distinct families of Laurent series
solutions :
$$
m_{1}\left( t\right)=\left\{\begin{array}{rl}
&\frac{\alpha_{1}}{t}+i\left( \alpha _{1}^{2}-2\right) \alpha
_{2}+ \circ\left( t\right),\\
&\frac{\alpha _{1}}{t}-i\left( \alpha _{1}^{2}-2\right) \alpha
_{2}+ \circ\left( t\right),
\end{array}\right.
$$
$$ m_{2}\left( t\right)= \left\{\begin{array}{rl} &\frac{i\alpha
_{1}}{t}-\alpha _{1}^{2}\alpha _{2}+\circ\left( t\right) ,\\
&\frac{-i\alpha _{1}}{t}-\alpha _{1}^{2}\alpha _{2} +\circ\left(
t\right) ,\end{array}\right.
$$
$$m_{3}\left( t\right)= \left\{\begin{array}{rl}
&\frac{i}{t}+\alpha _{1}\alpha _{2}
+\circ\left( t\right) ,\\
&\frac{-i}{t}+\alpha _{1}\alpha _{2}+\circ\left( t\right) ,
\end{array}\right.
$$
$$
\gamma _{1}\left( t\right)= \left\{\begin{array}{rl}
&\frac{1}{2t^{2}}+\circ\left( t\right) ,\\
&\frac{1}{2t^{2}}+\circ\left( t\right) ,\end{array}\right.
$$
$$\gamma _{2}\left( t\right)= \left\{\begin{array}{rl}
&\frac{i}{2t^{2}}+\circ\left( t\right) ,\\
&\frac{-i}{2t^{2}}+\circ\left( t\right) ,\end{array}\right.
$$
$$
\gamma _{3}\left( t\right) = \left\{\begin{array}{rl}
&\frac{\alpha _{2}}{t}+\circ\left( t\right) ,\\
&\frac{\alpha _{2}}{t}+\circ\left( t\right) ,\end{array}\right.$$
which depend on $5$ free parameters $\alpha _{1},...,$ $\alpha
_{5.}$ By substituting these series in the constants of the motion
$H_{i}$ (22), one eliminates three parameters linearly, leading to
algebraic relation between the two remaining parameters, which is
nothing but the equation of the divisor $\mathcal{D}$ along which
the $m_{i},\gamma _{i}$ blow up. Since the system (21) admits two
families of Laurent solutions, then $\mathcal{D}$ is a set of two
isomorphic curves of genus $3,$
$\mathcal{D}=\mathcal{D}_{1}+\mathcal{D}_{-1}:$
\begin{equation}\label{eqn:euler}
\mathcal{D}_{\varepsilon }:\text{ }P\left( \alpha _{1},\alpha
_{2}\right) =\left( \alpha _{1}^{2}-1\right) \left( \left( \alpha
_{1}^{2}-1\right) \alpha _{2}^{2}-P\left( \alpha _{2}\right)
\right)+c_{4}=0,
\end{equation}
where $$P\left( \alpha _{2}\right) =c_{1}\alpha
_{2}^{2}-2\varepsilon c_{2}\alpha _{2}-1,$$ and $\varepsilon =\pm
1.$ Each of the curve $\mathcal{D}_{\varepsilon }$ is a $2-1$
ramified cover $\left( \alpha _{1},\alpha _{2},\beta \right) $ of
elliptic curves $\mathcal{D}_{\varepsilon }^{0}:$
\begin{equation}\label{eqn:euler}
\mathcal{D}_{\varepsilon }^{0}:\beta ^{2}=P^{2}\left( \alpha
_{2}\right) -4c_{4}\alpha _{2}^{4},
\end{equation}
ramified at the $4$ points $\alpha _{1}=0$ covering the $4$ roots
of $P\left( \alpha _{2}\right) =0.$ It was shown by the author
$\left[14\right] $ that each divisor $\mathcal{D}_{\varepsilon }$
is ample and defines a polarization $\left( 1,2\right) ,$ whereas
the divisor $D,$ of geometric genus $9$, is very ample and defines
a polarization $\left( 2,4\right),$ accordind to (8). The affine
surface $$M_{c}=\bigcap_{i=1}^{4}\left\{ H_{i}=c_{i}\right\}
\subset \mathbb{C}^{6},$$ defined by putting the four invariants
(22) of the Kowalewski flow (21) equal to generic constants, is
the affine part of an abelian surface $\widetilde{M_{c}}$ with
\begin{eqnarray}
\widetilde{M_{c}}\text{ }\backslash \text{
}M_{c}=\mathcal{D}&=&\text{one
genus 9 curve consisting of two genus 3 }\nonumber\\
&&\text{curves }\mathcal{D}_{\varepsilon }\text{ }(23) \text{
intersecting in 4 points. Each }\nonumber\\
&&\mathcal{D}_{\varepsilon }\text{ is a double cover of an
elliptic curve
}\mathcal{D}_{\varepsilon }^{0}\text{ }(24) \nonumber\\
&&\text{ramified at 4 points.}\nonumber
\end{eqnarray}
Moreover, the Hamiltonian flows generated by the vector fields
$X_{H_{1}}$ and $X_{H_{4}}$ are straight lines on
$\widetilde{M_{c}}.$ The 8 functions $$1,\quad f_{1}=m_{1},\quad
f_{2}=m_{2},\quad f_{3}=m_{3},\quad f_{4}=\gamma _{3},\quad
f_{5}=f_{1}^{2}+f_{2}^{2},$$$$ f_{6}=4f_{1}f_{4}-f_{3}f_{5},\quad
f_{7}=\left( f_{2}\gamma _{1}-f_{1}\gamma _{2}\right)
f_{3}+2f_{4}\gamma _{2},$$ form a basis of the vector space of
meromorphic functions on $\widetilde{M_{c}}$ with at worst a
simple pole along $\mathcal{D}$ Moreover, the map
$$\widetilde{M_{c}}\simeq \mathbb{C}^{2}/Lattice\rightarrow
\mathbb{CP}^{7}\text{ },\text{ }\left( t_{1},t_{2}\right) \mapsto
\left[ \left( 1,f_{1}\left( t_{1},t_{2}\right) ,...,f_{7}\left(
t_{1},t_{2}\right) \right) \right] ,$$ is an embedding of
$\widetilde{M_{c}}$ into $\mathbb{CP}^{7}.$ Following the method
(subsection 3.1), we obtain the following theorem :
\begin{Theo}
The tori $\widetilde{M_{c}}$ can be identified as
$$\widetilde{M_{c}}=Prym^*(\mathcal{D}_\varepsilon/\mathcal{D}_\varepsilon
^0),$$ i.e., dual of $Prym(\mathcal{D}_{\varepsilon
}/\mathcal{D}_{\varepsilon }^{0})$ and the problem linearizes on
this Prym variety.
\end{Theo}

\subsection{Kirchhoff's equations of motion of a solid in an ideal fluid}

The Kirchhoff's equations of motion of a solid in an ideal fluid
have the form
\begin{eqnarray}
\dot p_{1}&=&p_{2}\frac{\partial H}{\partial l_{3}}-p_{3}
\frac{\partial H}{\partial l_{2}},\nonumber\\
\dot p_{2}&=&p_{3}\frac{\partial H}{\partial l_{1}}-p_{1}
\frac{\partial H}{\partial l_{3}},\nonumber\\
\dot p_{3}&=&p_{1}\frac{\partial H}{\partial l_{2}}-p_{2}
\frac{\partial H}{\partial l_{1}},\\
\dot l_{1}&=& p_{2}\frac{\partial H}{\partial
p_{3}}-p_{3}\frac{\partial H}{\partial p_{2}}+l_{2}\frac{\partial
H}{\partial
l_{3}}-l_{3}\frac{\partial H}{\partial l_{2}},\nonumber\\
\dot l_{2}&=&p_{3}\frac{\partial H}{\partial
p_{1}}-p_{1}\frac{\partial H}{\partial p_{3}}+ l_{3}\frac{\partial
H}{\partial
l_{1}}-l_{1}\frac{\partial H}{\partial l_{3}},\nonumber\\
\dot l_{3}&=& p_{1}\frac{\partial H}{\partial
p_{2}}-p_{2}\frac{\partial H}{\partial p_{1}}+ l_{1}\frac{\partial
H}{\partial l_{2}}-l_{2}\frac{\partial H}{\partial
l_{1}},\nonumber
\end{eqnarray}
where $(p_{1},p_{2},p_{3})$\ is the velocity of a point fixed
relatively to the solid, $(l_{1},l_{2},l_{3})$\ the angular
velocity of the body expressed with regard to a frame of reference
also fixed relatively to the solid and $H$\ is the hamiltonian.
These equations \ can be regarded as the equations of the
geodesics of the right-invariant metric on the group $E\left(
3\right) =SO\left( 3\right) \times \mathbb{R}^{3}$\ of motions of
3-dimensional euclidean space $\mathbb{R}^{3},$\ generated by
rotations and translations. Hence the motion has the trivial
coadjoint orbit invariants $\langle p,p\rangle$ and $\langle
p,l\rangle.$ As it turns out, this is a special case of a more
general system of equations written as
\begin{eqnarray}
\dot x&=&x\wedge \frac{\partial H}{\partial x}+y\wedge
\frac{\partial H}{\partial y},\nonumber\\
\dot y&=&y\wedge \frac{\partial H}{\partial x} +x\wedge
\frac{\partial H}{\partial y},\nonumber
\end{eqnarray}
where $x=\left( x_{1},x_{2},x_{3}\right) \in \mathbb{R}^{3}$ et
$y=\left( y_{1},y_{2},y_{3}\right) \in \mathbb{R}^{3}.$ The first
set can be obtained from the second by putting
$(x,y)=(l,p/\varepsilon)$ and letting $\varepsilon\rightarrow 0.$
The latter set of equations is the geodesic flow on $SO(4)$ for a
left invariant metric defined by the quadratic form $H.$ In
Clebsch's case, equations (25) have the four invariants :
\begin{eqnarray}
H_{1}&=&H=\frac{1}{2}\left(
a_{1}p_{1}^{2}+a_{2}p_{2}^{2}+a_{3}p_{3}^{2}+b_{1}l_{1}^{2}+
b_{2}l_{2}^{2}+b_{3}l_{3}^{2}\right),\nonumber\\
H_{2}&=&p_{1}^{2}+p_{2}^{2}+p_{3}^{2},\\
H_{3}&=&p_{1}l_{1}+p_{2}l_{2}+p_{3}l_{3},\nonumber\\
H_{4}&=&\frac{1}{2}\left(
b_{1}p_{1}^{2}+b_{2}p_{2}^{2}+b_{3}p_{3}^{2}+\varrho \left(
l_{1}^{2}+l_{2}^{2}+l_{3}^{2}\right) \right) ,\nonumber
\end{eqnarray}
with $$\frac{a_{2}-a_{3}}{b_{1}}+\frac{a_{3}-a_{1}}{b_{2}}+
\frac{a_{1}-a_{2}}{b_{3}}=0,$$ and the constant $\varrho $
satisfies the conditions $$\varrho =\frac{b_{1}\left(
b_{2}-b_{3}\right) }{a_{2}-a_{3}}=\frac{b_{2}\left(
b_{3}-b_{1}\right) }{a_{3}-a_{1}}=\frac{b_{3}\left(
b_{1}-b_{2}\right) }{a_{1}-a_{2}}.$$ The system (25) can be
written in the form (11) with $m=6$; to be precise
\begin{equation}\label{eqn:euler}
\dot{x}=f\left( x\right) \equiv J\frac{\partial H}{\partial
x}\text{ },\text{
}x=(p_{1},p_{2},p_{3},l_{1},l_{2},l_{3})^{\intercal },
\end{equation}
where
$$
J=\left(\begin{array}{cc}
O&P\\
P&L
\end{array}\right),
$$
$$
P=\left(\begin{array}{ccc}
0&-p_{3}&p_{2}\\
p_{3}&0&-p_{1}\\
-p_{2}&p_{1}&0
\end{array}\right),
$$
$$
L=\left(\begin{array}{ccc}
0&-l_{3}&l_{2}\\
l_{3}&0&-l_{1}\\
-l_{2}&p_{1}&0
\end{array}\right).
$$
Consider points at infinity which are limit points of trajectories
of the flow. In fact, there is a Laurent decomposition of such
asymptotic solutions,
\begin{equation}\label{eqn:euler}
x\left( t\right) =t^{-1}\left( x^{\left( 0\right) }+ x^{\left(
1\right) }t+x^{\left( 2\right) }t^{2}+...\right),
\end{equation}
which depend on $\dim (phase$\ $space)-1=5$\ free parameters$.$\
Putting (28) into (27), solving inductively for the $x^{\left(
k\right) },$\ one finds at the $0^{th}$\ step a non-linear
equation, $$x^{\left( 0\right) }+f(x^{\left( 0\right)})=0,$$ and
at the $k^{th}$\ step, a linear system of equations,
$$
(L-kI) x^{(k)}= \left\{
\begin{array}{rl}
0 & \mbox{for}\quad k=1\\
\text{quadratic polynomial in }x^{\left( 1\right) },...,x^{\left(
k\right) }& \mbox{for}\quad k \geq 1,
\end{array}
\right.
$$
where $L$ denotes the jacobian map of the non-linear equation
above.\ One parameter appear at the $0^{th}$ step, i.e., in the
resolution of the non-linear equation\ and the $4$\ remaining ones
at the $k^{th}$\ step, $k=1,...,4.$ Taking into account only
solutions trajectories lying on the invariant surface
$$M_{c}=\overset{4}{\underset{i=1}{\bigcap }}\left\{ H_{i}\left(
x\right) =c_{i}\right\} \subset \mathbb{C}^{6},$$ we obtain
one-parameter families which are parameterized by a curve. To be
precise we search for the set of Laurent solutions which remain
confined to a fixed affine invariant surface, related to specific
values of $c_{1},c_{2},c_{3},c_{4},$ i.e.,
\begin{eqnarray}
\mathcal{D}&=&\bigcap_{i=1}^{4}\left\{ t^{0}-\text{coefficient of
}H_{i}\left( x\left( t\right) \right) =c_{i}\right\},\nonumber\\
&=&\text{an algebraic curve defined by}\nonumber\\
&&\theta ^{2}+c_{1}\beta ^{2}\gamma ^{2}+c_{2}\alpha ^{2}\gamma
^{2}+c_{3}\alpha ^{2}\beta ^{2}+c_{4}\alpha \beta \gamma =0,
\end{eqnarray}
where $\theta $ is an arbitrary parameter and where $$\alpha
=x_{4}^{\left( 0\right) },\beta =x_{5}^{\left( 0\right) },\gamma
=x_{6}^{\left( 0\right) },$$ parameterizes the elliptic curve
\begin{equation}\label{eqn:euler}
\mathcal{E}:\beta ^{2}=d_{1}^{2}\alpha ^{2}-1,\text{ }\gamma ^{2}=
d_{2}^{2}\alpha^{2}+1,
\end{equation}
with $d_{1},d_{2}$\ such that: $$d_{1}^{2}+d_{2}^{2}+1=0.$$ The
curve $\mathcal{D}$ is a $2$-sheeted ramified covering of the
elliptic curve $\mathcal{E}$. The branch points are defined by the
$16$\ zeroes of $c_{1}\beta ^{2}\gamma ^{2}+c_{2}\alpha ^{2}\gamma
^{2}+c_{3}\alpha ^{2}\beta ^{2}+c_{4}\alpha \beta \gamma $\ on
$\mathcal{E}$. The curve $\mathcal{D}$ is unramified at infinity
and by Hurwitz's formula, the genus of $\mathcal{D}$ is $9$. Upon
putting $\zeta \equiv \alpha ^{2}$, the curve $\mathcal{D}$ can
also be seen as a $4-$sheeted unramified covering of the following
curve of genus $3:$
$$C:\left( \theta ^{2}+c_{1}\beta ^{2}\gamma ^{2}+
\left( c_{2}\gamma ^{2}+c_{3}\beta ^{2}\right) \zeta \right) ^{2}-
c_{4}^{2}\zeta \beta ^{2}\gamma^{2}=0.$$ Moreover, the map $$\tau
:C\longrightarrow C,\quad(\theta ,\zeta )\longmapsto (-\theta
,\zeta ),$$ is an involution on $C$ and the quotient
$C_{0}=C/\tau$ is an elliptic curve defined by
$$C_{0}:\eta ^{2}=c_{4}^{2}\zeta \left( d_{1}^{2}d_{2}^{2}\zeta ^{2}+
\left( d_{1}^{2}-d_{2}^{2}\right) \zeta -1\right).$$ The curve $C$
is a double ramified covering of $C_{0},$ $$C\longrightarrow
C_{0},\quad(\theta ,\eta ,\zeta )\longmapsto (\eta ,\zeta ),$$
$$
C:\left\{\begin{array}{rl} &\theta ^{2}=-c_{1}\beta ^{2}\gamma
^{2}-\left( c_{2}\gamma ^{2}+c_{3}\beta ^{2}\right) \zeta +\eta \\
&\eta ^{2}=c_{4}^{2}\zeta \left( d_{1}^{2}d_{2}^{2}\zeta
^{2}+\left( d_{1}^{2}-d_{2}^{2}\right) \zeta -1\right).
\end{array}\right.
$$
Let $(a_{1},a_{2},a_{3},b_{1},b_{2},b_{3})$\ be a canonical
homology basis of $C$\ such that $$\tau \left( a_{1}\right)
=a_{3},\quad\tau \left( b_{1}\right)=b_{3},\quad\tau \left(
a_{2}\right)=-a_{2},\quad\tau \left( b_{2}\right)=-b_{2},$$ for
the involution $\tau .$ Using the Poincar\'{e} residu map [9], we
show that
$$\omega _{0}=\frac{d\zeta }{\eta },\text{ }\omega _{1}=
\frac{\zeta d\zeta }{\theta \eta },\text{ }\omega _{2}=
\frac{d\zeta }{\theta \eta},$$ form a basis of holomorphic
differentials on $C$\ and $$\tau ^{*}\left( \omega _{0}\right)
=\omega _{0},\quad\tau ^{*}\left( \omega _{k}\right)=-\omega
_{k},\quad\left( k=1,2\right).$$ Haine [10] shows that the flow
evolues on an abelian surface $\widetilde{M}_c\subseteq
\mathbb{CP}^7$ of period matrix $\left(
\begin{array}{llll}
2 & 0 & a & c \\
0 & 4 & c & b
\end{array}
\right),\text{ Im}\left(
\begin{array}{ll}
a & c \\
c & b
\end{array}
\right) >0$ and also identified $\widetilde{M}_c$ as Prym variety
$Prym(C/C_0)$.
\begin{Theo}
The abelian surface $\widetilde{M}_c$ can be identified as
$Prym(C/C_0)$. More precisely
$$
\overset{4}{\underset{i=1}{\bigcap }}\left\{x\in \mathbb{C}^{6},
H_{i}\left( x\right) =c_{i}\right\}=Prym(C/C_0) \backslash
\mathcal{D},$$ where $\mathcal{D}$ is a genus 9 curve (29), which
is a ramified cover of an elliptic curve $\mathcal{E}$ (30) with
16 branch points.
\end{Theo}

\end{document}